\documentclass[10pt]{article}

\usepackage[utf8]{inputenc}
\usepackage{amsmath} 
\usepackage{amssymb} 
\usepackage{amsthm}
\usepackage{xcolor}
\usepackage{parskip}
\usepackage[T1]{fontenc}
\usepackage[french]{babel}

\newtheorem{proposition}{Proposition}[subsection]
\newtheorem{Remarque}[proposition]{Remarque}

\newtheorem{fait}[proposition]{Fait}
\newtheorem{definition}[proposition]{Définition}
\newtheorem{theoreme}[proposition]{Théorème}
\newtheorem{Lemme}[proposition]{Lemme}
\newtheorem{corollaire}[proposition]{Corollaire}
\def\Ind#1#2{#1\setbox0=\hbox{$#1x$}\kern\wd0\hbox to
0pt{\hss$#1\mid$\hss}
\lower.9\ht0\hbox to 0pt{\hss$#1\smile$\hss}\kern\wd0}

\def\Notind#1#2{#1\setbox0=\hbox{$#1x$}\kern\wd0\hbox to 0pt{\mathchardef
\nn="3236\hss$#1\nn$\kern1.4\wd0\hss}\hbox to 0pt{\hss$#1\mid$\hss}\lower.9\ht0
\hbox to 0pt{\hss$#1\smile$\hss}\kern\wd0}
\def\nind{\mathop{\mathpalette\Notind{}}}

\title{Symétrons et K-boucles $\omega$-stables\footnote{ mots-clé : espace symétrique, K-boucle, rang de Morley, rang de Lascar
MSC 2020 : 20N05, 20F11, 03C60}}

\author{Samuel Zamour \footnote{Univ Paris Est Creteil, Univ Gustave Eiffel, CNRS, LAMA UMR8050, F-94010 Creteil, France, e-mail : samuel.zamour@u-pec.fr}}
\date{19 septembre 2023}
\begin{document}
\maketitle
\begin{abstract}
    Nous développons la théorie des modèles des K-boucles et des symétrons $\omega$-stables. En poursuivant le travail séminal de Poizat, nous établissons notamment une version appropriée du théorème des indécomposables et nous adaptons l'analyse de Lascar à ce contexte.
    
    \begin{center}
        \textbf{Abstract}
    \end{center}
  We develop the model theory of $\omega$-stable K-loops and 'symétrons'. Continuing Poizat's seminal work, we notably establish an appropriate version of the indecomposability theorem and we adapt Lascar's analysis to this context.

\end{abstract}

\section{Introduction}
Les espaces symétriques apparaissent naturellement dans la classification des groupes de Lie réels. Ce sont des variétés riemaniennes telles qu'en chaque point la réflexion géodésique soit une isométrie. Suivant les travaux pionniers de Loos \cite{Lo}, on peut en donner une axiomatisation algébrique en utilisant une opération binaire (voir \cite{Lo} p.63). On dit que $(S,s)$ est un \textit{espace symétrique} (différentiel), où $S$ est une variété différentielle et $s:S\times S \rightarrow S$ est une application différentiable, si :
\begin{enumerate}
    \item $s(x,x)=x.$
    \item $s(s(x,y),y)=x.$
    \item $s(s(x,y),z)=s(s(x,z),s(y,z)).$
    \item Pour tout $x$ dans $S$, il existe un voisinage $U$ de $x$ tel que $s(x,y)=y$ implique que $y=x$ pour tout $y$ dans $U$.
\end{enumerate}
Par exemple, un groupe de Lie muni de l'opération $s(x,y)=xy^{-1}x$ est un espace symétrique. On pourra consulter \cite{Lo} pour plus d'exemples.
\\
De façon indépendante, Bruno Poizat a réintroduit la notion d'espace symétrique (abstrait) afin d'analyser la classification des groupes de rang de Morley fini de rang 3, obtenue par O. Frécon \cite{Fre}. En particulier, il remarque qu'un groupe sans involutions de rang de Morley fini peut être muni d'une structure de \textit{symétron} avec l'opération $s(x,y)=xy^{-1}x$ (on considère un ensemble de base abstrait sans hypothèse de nature géométrique et on remplace la condition 4 par la condition  :
\[ \text{4'. Pour tous x,y dans S, il existe un unique z tel que}~y=s(x,z).)\]
 Les articles \cite{Poi3}, \cite{Poi4} et \cite{Poi5} jettent les bases d'une analyse des espaces symétriques dans la perspective de la théorie des modèles. L'$omega$-stabilité constitue un bon cadre pour mener à bien cette étude. Pour une référence générale sur la théorie des modèles, on peut consulter \cite{Las} ou \cite{TZ}.
\\
\\
Les espaces symétrique sont liés à des structures algébriques binaires appelées \textit{K-boucles} : il y a un élément neutre mais l'associativité de l'opération est seulement partielle. Ces structures dans le cas fini ont été étudiées en détail par Glauberman \cite{Glau}. Les K-boucles sont l'objet d'un regain récent d'attention en théorie des modèles. En particulier, les résultats de Tent et Clausen \cite{CT} sur les groupes de Frobenius de rang de Morley fini (y compris sans involutions) utilisent de façon cruciale la notion de K-boucle. 
\\
Nous proposons dans cet article d'exploiter le lien entre espaces symétriques et K-boucles afin de prolonger les résultats de Poizat et d'établir une analyse modèle-théorique commune de ces structures.
La première section rappelle certaines définitions et propriétés concernant les boucles. Nous nous intéressons d'abord aux conditions de chaîne et aux K-boucles minimales (section 2). Sur la base de \cite{Nagy}, nous établissons ensuite la co-interprétabilité entre les symétrons et les K-boucles uniquement 2-divisibles (section 3). Dans ce cadre, nous donnons les éléments d'une théorie des modèles commune aux symétrons et aux K-boucles $\omega$-stables en établissant en particulier une nouvelle forme du théorème des indécomposables (section 4). Enfin, dans une dernière section (section 5), nous adaptons l'analyse de Lascar à notre contexte.
\\
Cet article se base sur un chapitre de la thèse doctorale de l'auteur \cite{Zam}, effectuée sous la direction de Frank Wagner au sein de l'Université Claude Bernard Lyon I. Nous souhaiterions le remercier chaleureusement pour ses relectures attentives et ses nombreuses remarques qui se sont avérées essentielles en bien des points.
Nous remercions également Bruno Poizat pour ses encouragements; nous lui sommes reconnaissants de nous avoir aimablement communiqué les versions préliminaires de certains de ses articles.
\section{Préliminaires sur les boucles}
Débutons par quelques rappels concernant les boucles. Pour une présentation exhaustive des notions pertinentes, nous renvoyons à \cite{Kie} et \cite{Bru}.
\begin{definition}
    Soit $(B,+,0)$ une structure munie d'une opération binaire $+$. 
    On dit que $B$ est une boucle si :
    \begin{itemize}
        \item Pour tout a dans $B$, on a : \[0+a=a+0=a.\]
        \item Les équations $x+a=b$ et $a+y=b$ ont chacune une unique solution dans $B$ pour tous $a,b$ dans $B$.
    \end{itemize}
    Si $B$ satisfait de plus les conditions suivantes, on dit que $B$ est une K-boucle :
    \item (Condition de Bol) Pour tous $a,b,c\in B$, on a : \[a+(b+(a+c))=(a+(b+a))+c.\]
\item (Propriété de l'inversion par automorphisme) Pour tous $a,b\in B$, on a : \[-(a+b)=-a-b.\]
\end{definition}
\begin{Remarque}
Si une boucle satisfait seulement la condition de Bol, on parle de boucle de Bol. Notons que cette condition est suffisante pour démontrer l'existence d'un unique inverse.
\end{Remarque}
On dit qu'une K-boucle est \textit{uniquement 2-divisible} si l'application $x\mapsto x+x$ définit une bijection. Voici quelques exemples de K-boucles (uniquement 2-divisibles) :
\begin{enumerate}
    \item Un groupe uniquement 2-divisible avec l'opération $x+y=x^{1/2}yx^{1/2}$.
    \item Si $\alpha$ est un automorphisme involutif d'un groupe $G$, l'ensemble $X=\{ \alpha(g)g^{-1} : g\in G\}$ peut être muni d'une structure de K-boucle si $X$ est uniquement 2-divisible.
    \end{enumerate}
Pour développer la théorie des boucles, il faut disposer d'une notion appropriée de sous-boucle et de boucle quotient. Contrairement à ce qui passe dans le cas des groupes, définir un concept opératoire de sous-structure normale et plus généralement d'automorphisme intérieur pose problème.
\\
Nous commençons par présenter une première notion d'application intérieure.
Soit $(B,+,0)$ une boucle, et soit $S_B$ son groupe de permutation. Pour tout $a$ dans $B$, on définit la \textit{translation à gauche} $g_a : x\rightarrow a+x$ et on considère le groupe engendré par les translations (à gauche) : \[ M_g(B)=\langle g_a : a\in B \rangle\leq S_B. \] 
Remarquons que les translations sont des bijections définissables (et donc leurs inverses sont aussi des bijections définissables).  Nous pouvons dès lors définir les \textit{applications de précession} (à gauche)  $\delta_{a,b}=g_{a+b}^{-1}g_ag_b$. Elles sont caractérisées par la propriété suivante  :
\[ a+(b+c)=(a+b)+\delta_{a,b}(c). \]
\\
Le \textit{groupe des applications intérieures} (à gauche) est alors :
\[D(B)=\langle \delta_{a,b} : a,b\in B \rangle\leq S_B. \]
\begin{Remarque}
On peut aussi définir la translation à droite $r_b : x\longrightarrow x+b$ pour $b\in B$ et ainsi considérer le groupe des translations bilatères :
\[ M(B)=\langle g_a,r_b : a,b\in B\rangle. \]
\end{Remarque}
Dans les K-boucles, un certain nombre de propriétés sont satisfaites, ce qui facilite les calculs. En particulier, l'opération $a\cdot n =\underset{n fois}{\underline{a+...+a}}$ est bien définie pour tout $n\in \mathbb{N}$ (et donc pour tout $n\in\mathbb{Z}$).
\begin{fait} (voir \cite{Kie})\label{boucle général}
\begin{enumerate}
    \item $a+b=\delta_{a,b}(b+a)$.
    \item Pour tout $m,n\in \mathbb{Z}$, $a\cdot n+(a\cdot m+x)=a\cdot (n+m)+x$.
    \item Pour tout $m,n\in \mathbb{Z}$, $\delta_{a\cdot m,a\cdot m}=Id$.
    \item $\delta_{a,b}^{-1}=\delta_{b,a}$.
    \item $\delta_{a,b}=\delta_{-b,b+a}$.
    \item $\delta_{a,b+a}=\delta_{a,b}$.
    \item Les applications de précession sont des automorphismes de la K-boucle.
    \item $(a+b)\cdot 2=a+(b\cdot 2+a)$
    \item L'application $x\mapsto x\cdot 2$ est injective ssi $K$ ne contient pas d'élément d'ordre 2.
    \item Si $B$ est une K-boucle uniquement 2-divisible et si $\epsilon$ est un automorphisme involutif sans point fixe de $B$ alors $\epsilon(x)=-x$ pour tout $x\in B$. 
\end{enumerate}
\end{fait}
\begin{Remarque} Les trois premières propriétés restent vraies lorsque l'on considère seulement une boucle de Bol. En particulier, tout élément d'une boucle de Bol engendre un groupe cyclique.
\end{Remarque}
Dans les boucles de Bol (et donc \textit{a fortiori} dans les K-boucles), les sous-boucles peuvent être caractérisées assez aisément :
\begin{fait}\cite[Lemme 1.3]{Nagy}\label{sous-boucle}
Soit $B$ une boucle de Bol. Alors $C\subseteq B$ est une sous-boucle ssi $C$ est stable par inversion et addition.
\end{fait}
\begin{Remarque} Les sous-boucles d'une K-boucle sont des K-boucles.
\end{Remarque}
En général, la notion de quotient d'une K-boucle $B$ par une sous K-boucle $C$ n'a pas de sens : la relation $x\equiv y$ ssi $x-y\in C$ n'est pas forcément une relation d'équivalence. Néanmoins, sous certaines conditions, un quotient bien défini peut être muni d'une structure de boucle. D'après le paragraphe IV.1 de \cite{Bru}, une sous-boucle $C$ est \textit{normale} si elle satisfait les conditions suivantes :
\[ b+C=C+b,\quad (a+b)+C=a+(b+C),\quad C+(a+b)=(C+a)+b \]
On peut également caractériser les sous-boucles normales en introduisant une notion plus forte de groupe d'applications intérieures.
\begin{definition}\label{groupe application intérieure}
Soit $B$ une boucle. On définit $I(B)$, le groupe des applications intérieures, comme le sous-groupe de $M(B)=\langle g_a, r_b : a,b\in B \rangle$, engendré par les applications suivantes :
\[ r_{a,b}=r_ar_br_{a+b}^{-1},\quad g_{a,b}=g_ag_bg_{b+a}^{-1},\quad c_a=r_a g_{a}^{-1}.\]
\end{definition}
\begin{fait}\cite[Sec.1. Chap IV]{Bru} Soit $B$ une boucle.
\begin{enumerate}
    \item Le groupe $I(B)$ correspond exactement à l'ensemble des applications $\alpha\in M(B)$ qui fixe $0$.
    \item Une sous-bouble $C\leq B$ est normale ssi elle est stable sous l'action de $I(B)$.
\end{enumerate} 
\end{fait}
Le noyau d'un homomorphisme de boucles est toujours une sous-boucle normale. Les résultats valables pour les groupes se généralisent alors dans notre contexte. En particulier, le quotient d'une boucle par une sous-boucle normale est bien défini et on peut le munir d'une structure de boucle compatible : 
\begin{fait}\cite[Section IV.1]{Bru}
\begin{enumerate}
    \item Si $C$ est une sous-boucle normale de $B$, alors l'application $x\mapsto x+C$ définit un homomorphisme de boucles entre $B$ et $B/C$. De plus, si $C$ est le noyau d'un homomorphisme de boucles $\phi$, alors $\phi$ induit un isomorphisme naturel entre $B/C$ et $\phi(B)$.
    \item Si $C$ est une sous-boucle engendrée par un ensemble de sous-boucles normales, alors $C$ est une sous-boucle normale.
    \item Si $C$ et $D$ sont des sous-boucles telles que $D$ est normale dans la sous-boucle $\langle D,C \rangle$ engendrée par $D$ et $C$, alors $\langle D,C \rangle=C+D=D+C$, $C\cap D$ est une sous-boucle normale de $C$; $C+D/D\simeq C/(C\cap D)$.
\end{enumerate}
\end{fait}
Nous introduisons une dernière notion importante dans la théorie des K-boucles.
\begin{definition}\label{def sans point fixe}
On dit qu'une K-boucle $B$ est sans point fixe, si toute application intérieure (à gauche) non-trivial $\phi\in D(B)$ n'a pas de point fixe non-trivial.
\end{definition}
\begin{Remarque}
Les K-boucles qui apparaissent dans la classification des groupes strictement 2-transitifs sont sans point fixe (voir \cite{KW 1}, \cite{KW 2} et \cite{KW 3}).
\end{Remarque}
Dans une K-boucle sans point fixe $B$, un automorphisme de précession $\delta_{a,b}$ est définissable sur $(a+b),(b+a)$.
\\
En effet, si $\delta_{x,y}(b+a)=\delta_{a,b}(b+a)=(a+b)$ alors $\delta_{b,a}\delta_{x,y}(b+a)=b+a$ et donc $\delta_{b,a}\delta_{x,y}$ admet un point fixe. Puisque $\delta_{b,a}\delta_{x,y}\in D(B)$, on obtient que $\delta_{a,b}=\delta_{x,y}$.
\begin{Remarque}
Dans cet article, nous n'étudions pas la notion de demi-plongement qui constitue la clef de voûte de l'analyse de Glauberman. Une K-boucle $B$ est demi-plongée dans un groupe $G$ s'il existe une injection $\eta : B\longrightarrow G$ telle que $\eta(B)\subseteq G$ satisfait les conditions suivantes :
\[ 1\in \eta(B),\quad \eta(B)^{-1}\subseteq \eta(B),\quad a\eta(B)a\subseteq \eta(B)\quad \mbox{pour tout}\quad a\in\eta(B). \]
On suppose de plus que l'application d'élévation au carré (au sens de l'opération du groupe) est bien définie et bijective sur $\eta(B)$.
Les K-boucles uniquement 2-divisibles peuvent être demi-plongées dans un groupe mais ce dernier n'est pas en général interprétable et nous ne pouvons pas en extraire des informations pertinentes du point de vue de la théorie des modèles.
\end{Remarque}

\section{Premières investigations modèle-théoriques et K-boucles minimales}
Quelles propriétés pertinentes du point de vue de la théorie des modèles peut-on établir pour les K-boucles $\omega$-stables ?
\\
L'étude des conditions de chaîne constitue un préalable indispensable. On note $RM$ le rang de Morley et $DM$ le degré de Morley.
\\
Soient $B$ une K-boucle $\omega$-stable et $C$ une sous-boucle définissable. Soit $x\in B\setminus C$, on remarque que $C\cap C+x=\varnothing$ :  si $c+x=c'$ alors $-c+(c+x)=x=-c+c'$. Notons de plus que $RM(C)=RM(C+x)$.
\\
Par conséquent, ou bien le rang de $C$ ou bien son degré est strictement inférieur. Nous venons donc d'établir une condition de chaîne descendante sur les sous-boucles définissables. En particulier, toute partie non-vide est contenue dans une plus petite sous-boucle définissable.
\\
Si une sous-boucle est un groupe abélien, son enveloppe définissable est-elle un groupe abélien? Nous allons voir que c'est bien le cas pour le groupe cyclique engendré par un élément si on suppose que la K-boucle est sans point fixe et $\omega$-stable. 
\\
Soit $B$ une K-boucle sans point fixe (Définition \ref{def sans point fixe}). Rappelons que dans ce cas $b+b'= b'+b$ est équivalent à $\delta_{b,b'}=\delta_{b',b}=Id$.
\\
On définit pour un élément $x\in B$, le \textit{centralisateur} : \[C_B(x)=\{ b\in B : \delta_{x,b}=\delta_{-x,b}=Id \}\] 

 Nous pouvons maintenant introduire le \textit{centre} du centralisateur : \[Z(C_B(x))=\{b\in C_B(x) : \delta_{b,b'}=Id,\quad{pour\ tout\ } b'\in C_B(x) \} \]
\begin{Lemme}
Soit $B$ une K-boucle sans point fixe. L'ensemble $C_B(x)$ est stable par l'application $y\mapsto y\cdot n$, pour tout $n\in \mathbb{Z}$. De plus, $Z(C_B(x))$ est un sous-groupe abélien définissable contenant le groupe engendré par $x$.
\end{Lemme}
\begin{proof}
Tout d'abord, si $b\in C_B(x)$, alors $b-x=-x+b$ et donc $-b+x=x-b$ (propriété de l'inversion par automorphisme).
On procède par récurrence. Supposons que $b\cdot n+x=x+b\cdot n$. On a :
\[ b\cdot (n+1)+x=(b+b\cdot n)+x=b+(b\cdot n+x)=b+(x+b\cdot n)\]
\[ =(b+x)+\delta_{b,x}(b\cdot n)=(b+x)+b\cdot n=(x+b)+b\cdot n=x+b\cdot (n+1). \]
Les autres cas se traitent de façon analogue. On remarque que pour tout $b\in C_B(x)$, on a : $b+x\cdot n=x\cdot n+b$ pour $n\in \mathbb{Z}$. Puisque $x$ engendre un groupe abélien, on a bien $\langle x \rangle\subseteq Z(C_B(x))$.
\\
\\
Il suffit de vérifier que $Z(C_B(x))$ est stable par addition et par inversion (Proposition \ref{sous-boucle}).
\\
Soit $b_1\in Z(C_B(x))$ et soit $b\in C_B(x)$, on note tout d'abord que $b_1-b=-b+b_1$ car $-b\in C_B(x)$. Par conséquent, $-b_1+b=b-b_1$ (propriété de l'inversion par automorphisme).
\\
Soient $b_1,b_2\in Z(C_B(x))$ et $b\in C_B(x)$. On doit montrer que $b_1+b_2$ commute avec $b$ (c'est suffisant car la boucle est sans point fixe). Notons que $\delta_{b,b_1}=\delta_{b,b_2}=\delta_{b_1,b_2}=Id$.
\[ b+(b_1+b_2)=(b+b_1)+b_2=(b_1+b)+b_2=b_1+(b+b_2)=b_1+(b_2+b)=(b_1+b_2)+b. \]
Il s'agit donc d'une sous-boucle définissable mais puisque tous ses éléments commutent entre eux, c'est un groupe commutatif.
\end{proof}

\begin{corollaire}
Soit $B$ une K-boucle sans point fixe qui contient un élément d'ordre infini. Alors $B$ contient un sous-groupe abélien définissable infini.
\end{corollaire}
\begin{proof}
Il suffit de considérer $Z(C_B(x))$ pour un élément $x$ d'ordre infini.
\end{proof}
\begin{corollaire}\label{relèvement torsion boucle}
Soit $B$ une K-boucle sans point fixe et $\omega$-stable. Alors pour tout $x\in B$, la sous-boucle $d(x)=d(\langle x \rangle)$ est un groupe commutatif. \end{corollaire}
Nous pouvons également caractériser les K-boucles minimales.
\begin{corollaire}
Soit $B$ une K-boucle sans point fixe  $\omega$-stable qui contient un élément d'ordre infini. Si $B$ est minimale, i.e., ses sous-boucles définissables propres sont finies, alors $B$ est un groupe abélien.
\end{corollaire}
\begin{Remarque}
Dans \cite{CT}, les deux auteurs montrent qu'une K-boucle uniquement 2-divisible fortement minimale demi-plongée dans un groupe rangé $G$ est un groupe abélien.
\end{Remarque}
\section{Co-interprétabilité entre K-boucle et symétron}
Précisons maintenant les liens entre les symétrons et les K-boucles. 
 \begin{definition}\cite{Poi4}
Soit $X$ un ensemble et $s(x,y)$ une opération binaire. On dit que $(X,s)$ est un espace symétrique si l'opération $s$ satisfait les conditions suivantes :
\begin{enumerate}
    \item $s(x,x)=x$
    \item $s(s(x,y),y)=x$
    \item $s(s(x,z), s(y,z))= s(s(x,y),z)$
    \end{enumerate}
\end{definition}
\begin{definition}
On notera $s_y(x)=s(x,y)$ et on parlera de symétrie de centre $y$. Un symétron est un espace symétrique tel que pour tous $x,y$ il existe un unique $z$ qui vérifie $s_z(x)=y$. On dira alors que $z=m(x,y)$ est le milieu entre $x$ et $y$. 
\end{definition}
Remarquons que les applications  $s(.,y)$ et $m(.,y)$ sont des bijections (définissables).
\begin{proposition} Soit $(B,+)$ une K-boucle uniquement 2-divisible. Alors c'est un symétron pour l'opération  $s(x,y)= y+(-x+y)$. 
\end{proposition}
\begin{proof}
On  vérifie successivement les quatres conditions.
\[s(x,x)=x+(-x+x)=x\]
De plus,
\[ \begin{aligned}
    s(s(x,y),y)
&=y+[-(y+(-x+y))+y]\\
&=y+[(-y+(x-y))+y]\\
&=y+[-y+(x+(-y+y)]\\
&=y+(-y+x)=x.
\end{aligned} \]
\\
Pour ce qui est de la troisième condition, on obtient : 
\[ \begin{aligned}s(s(x,z),s(y,z))
&=(z+(-y+z)+[-(z+(-x+z))+(z+(-y+z)]\\
&= (z+(-y+z))+[(-z+(x-z))+(z+(-y+z)]\\
&= (z+(-y+z))+[-z+(x+(-z+(z+(-y+z))]\\
&= (z+(-y+z))+[-z+(x+(-y+z))]\\
&= z+[-y+(z+(-z+(x+(-y+z))]\\
&= z+[-y+(x+(-y+z))]\\
&= z+[(-y+(x-y))+z]\\
&= s(s(x,y),z). \end{aligned} \]
\\
Montrons finalement qu'entre deux points $a$ et $b$, il existe un unique milieu. 
\\
Soit $a\in B$; on note que $s(0,a/2)=a/2+(-0+a/2)=a$ et donc $a/2$ est le milieu entre $0$ et $a$ (par unique 2-divisibilité). On veut trouver une unique solution $z$ à l'équation $s(a,z)=b$. Or, par unique 2-divisibilité, l'équation $s(0,z)=s(b,a/2)$ a une unique solution, $z=s(b,a/2)/2$. Montrons que $s(a,s(s(b,a/2)/2,a/2))=b$ : 
\[ \begin{aligned}s(a,s(s(b,a/2)/2,a/2))
&=s(s(0,a/2),s(s(b,a/2)/2),a/2))\\
&=s(s(0,s(b,a/2)/2),a/2))\\
&=s(s(b,a/2),a/2)=b.\end{aligned}\]\qedhere
\end{proof}
\begin{Remarque}  Dans \cite{Nagy},  les auteurs montrent qu'une boucle de Bol uniquement 2-divisible est déjà un espace symétrique (proposition 9.2). Nous reprenons la démonstration qu'ils donnent de l'existence d'un unique milieu.
\end{Remarque}
De façon réciproque, on peut définir une structure de K-boucle uniquement 2-divisible à partir de la structure de symétron. En fait, la notion de symétron au sens de \cite{Poi4} coïncide avec la notion de \textit{quasi-groupe symétrique} étudiée dans \cite[Sec. 9-11]{Nagy}. En effet, $(X,s)$ est un symétron ssi $(X,+)$ est un quasi-groupe symétrique où $x+y=s(y,x)$. Le fait suivant nous indique la façon d'interpréter la structure de K-boucle.
\begin{fait}\cite[Proposition 9.11 et Remarque 10.1]{Nagy}
Soit $a\in (X,+)$ où $(X,+)$ est un quasi-groupe symétrique. Notons $L_y(x)=y+x=s(x,y)$ et $R_y(x)=x+y=s(y,x)$ et considérons ces applications comme des éléments du groupe symétrique de $X$. Alors l'opération \[ x\underset{a}{+}y=R_a^{-1}(x)+L_a(y)=s(L_a(y),R_a^{-1}(x))=s(s(y,a),m(x,a)). \]
définit sur $X$ une structure de K-boucle uniquement 2-divisible (on note $(X,s^a)$ le symétron associé). De plus, l'application $s_u$ définit un isomorphisme de boucles entre $(X,\underset{a}{+})$ et $(X,\underset{s_u(a)}{+})$. 
\\
L'application $m(x,a)$ définit quant à elle un isomorphisme de symétrons entre $(X,s^a)$ et $(X,s)$. 
\end{fait}
Dans la mesure où la structure de symétron et celle de K-boucle (uniquement 2-divisible) sont co-interprétables, les propriétés caractéristiques de la théorie des modèles ---  stablité, RM fini, ... --- se transfèrent d'un point de vue à l'autre.
\\
En revanche, le groupe dans lequel une K-boucle uniquement 2-divisible est demi-plongée n'est pas forcément interprétable à partir de la structure de boucle (ni de symétron). 

\section{Théorème des indécomposables et symétriseurs}
Dans cette section, nous rappelons les principaux éléments de la théorie des modèles des symétrons développée par B. Poizat dans \cite{Poi4}. Comme contribution originale, nous généralisons la version du théorème des indécomposables donnée pour les symétrons de rang de Morley fini.
 
\subsection{Conditions de chaîne, symétriseurs et types génériques}
Commençons par deux définitions : 
\begin{definition}
Soit $Y$ un sous-ensemble du symétron $S$. Le sous-ensemble $Y$ est une partie convexe si $s_y(x)\in Y$, pour tous $x,y\in Y$. On dira que $Y$ est clos par prise de milieu si  pour tous $x,y\in Y$, on a $m(x,y)\in Y$. Enfin, le symétriseur de $Y$, $Sym(Y)$, est l'ensemble des points $x$ tel que $s_x(Y)=Y$.
\end{definition}
On remarque que le symétriseur est un ensemble convexe. Soit $x,x'\in Sym(Y)$, et $y\in Y$. Il existe donc $y'\in Y$ tel que $y=s_x(y')$. Par conséquent,
\[s(y,s(x',x))=s(s(y',x),s(x',x))=s(s(y',x'),x)\in Y\]
car $s(y',x')\in Y$. On peut inverser les rôles de $x$ et $x'$.
\\
De plus, $Sym(Y)$ s'injecte dans $Y$ si $Y$ est non-vide : pour un $y\in Y$ fixé, à chaque $x\in Sym(Y)$, on associe l'élément $s_x(y)$ qui appartient à $Y$.
\\
Pour le moment, nous n'avons pas fait d'hypothèses qui relèvent de la théorie des modèles. Néanmoins, dans le cas $\omega$-stable et \textit{a fortiori} si le rang de Morley est fini, la théorie des symétrons présente des aspects remarquables. 
\begin{fait}\cite[Théorème 4]{Poi4}\label{sous-symétron}
Soit $X$ une partie définissable du symétron $\omega$-stable $S$. Alors $X$ est convexe ssi $X$ est clos par prise de milieu. De plus, $X$ est égal à son symétriseur. On parlera de sous-symétron.
\end{fait}
Notons que la propriété de convexité est préservée par intersection finie. 
\begin{Remarque}
Il n'y a pas de chaîne infinie descendante de sous-symétrons définissables.
\end{Remarque}
\begin{proof}
En effet, soit $Y\subset S$ un sous-symétron définissable, et soit $x\in S\setminus Y$. Pour tout $y\in Y$, $s_y(x)\in S\setminus Y$ (car $Y$ est clos par les symétries ainsi définies et ce sont des applications involutives) et donc l'application qui à $y$ associe $s_y(x)$ définit une injection de $Y$ vers $S\setminus Y$. Par conséquent, ou bien $RM(Y)<RM(S)$ ou bien $RM(Y)=RM(S)$ mais $DM(Y)<DM(S)$.
\end{proof}
Dans le contexte d'un symétron $\omega$-stable, il est possible d'étudier les types, et en particulier les types génériques (\textit{i.e.} de rang maximal) du point de vue de leur symétriseur.
\\
Soient $p$ et $q$ deux types; on note $Sym(p\longleftrightarrow q)$ l'ensemble des $x$ tels que $s_x(p)=q$. En raison de la définissabilité des types et de la condition de chaîne descendante sur les sous-symétrons définissables, il s'agit d'un ensemble définissable.
\\
Les ensembles $Sym(p)=Sym(p\longleftrightarrow p)$ et $Sym(p_1,...,p_n)$ où les $p_i$ sont permutés par symétrie sont des convexes définissables. En effet, si $x$ et $y$ sont dans $Sym(p_1,...,p_n)$, alors pour une réalisation $a_i$ de $p_i$, \[s(a_i,s(x,y))=s(s(a_j,y),s(x,y))=s(s(a_j,x),y)=s(a_k,y)=a_m\] avec $a_j,a_k,a_m$ des réalisations de types appartenant à l'ensemble $\{p_1, ..., p_n \}$.
\\
\\
On remarque qu'étant donné deux génériques indépendants $x$ et $y$, le milieu $z$ est un générique indépendant de $x$ et de $y$, de même $x$ et $s_y(x)$ sont des génériques indépendants. Nous obtenons une caractérisation des génériques similaire à la situation rencontrée pour les groupes $\omega$-stables :
\begin{fait}\cite[Théorèmes 12 et 13]{Poi4}
\begin{enumerate}
    \item Dans un symétron $S$ $\omega$-stable, toute partie convexe définissable non-vide est le symétriseur de ses types génériques; elle se décompose de manière unique en un nombre fini, qui est impair, de sous-ensembles définissables convexes disjoints de degré de Morley égal à un.
    \item Une partie $X$ définissable est générique ($RM(X)=RM(S)$) ssi $S$ est recouvert par un nombre fini de translatés de $X$, \textit{i.e} de parties de la forme $s_us_vX$.  \end{enumerate}
\end{fait}
Nous pouvons obtenir un analogue du Lemme 2.3 de \cite{Poi1} dans le contexte des symétrons $\omega$-stables :
\begin{Lemme}\label{symétrique générique}
Soit $S$ un symétron $\omega$-stable. On considère $S'=Sym(p_1,...,p_n)$ (pour des types $p_i$ de même rang). Alors $RM(S')\leq RM(p_i)$, $1\leq i \leq n$. Si $S'=Sym(p)$ et si $RM(p)=RM(S')$, alors un générique de $S'$ est au mileu de $p$ et d'un élément $c\in S$.
\end{Lemme}
\begin{proof}
On remarque que $S'$ est le symétriseur de ses types génériques. On réalise $p_i$ par $a$ et un type générique de $S'$ par $b$, de manière indépendante. On a donc :
\[ RM(b/S)=RM(b/S\cup a)=RM(s_b(a)/S\cup a)\leq RM(s_b(a)/S)=RM(p_j)=RM(p_i), \]
où on obtient la deuxième égalité en remarquant que chaque point du triplet $\{a,b,s_b(a)\}$ est définissable sur les deux autres.
\\
\\
Supposons désormais que $S'=Sym(p)$. En cas d'égalité, le $tp(s_b(a)/S\cup a)$ est l'héritier de $tp(s_b(a)/S)$. Or, $m(s_b(a),a)=b\in Sym(p)$ donc par héritage il existe $c\in S$ tel que $m(s_b(a), c)\in Sym(p)$. Puisque $s_b(a)$ réalise $p$, le type "milieu" $m(p,c)$ satisfait la formule définissant $Sym(p)$ et donc il s'agit d'un générique de $Sym(p)$ (car il est de même rang que $p$).
\end{proof}
\subsection{Le théorème des indécomposables}
Dans les symétrons de RM fini, il est possible de démontrer un analogue du théorème des indécomposables. Nous généralisons le théorème obtenu par B. Poizat dans \cite{Poi4}.
\begin{definition}
On dit qu'un symétron $\omega$-stable est connexe s'il ne contient pas de sous-symétron définissable propre de même rang. 
\end{definition}
\begin{definition}
Soit $A$ une partie définissable d'un symétron. On dit qu'elle est indécomposable si pour toute famille finie d'ensembles convexes définissables deux-à-deux disjoints $(X_i : 1\leq i \leq n)$, on a : $A\subseteq X_1\cup...\cup X_n$ implique $A\subseteq X_i$ pour un certain $i\in \{1, ..., n \}$.
\end{definition}
Le lemme suivant montre que l'indécomposabilité est préservée par prise de symétrie :
\begin{Lemme}
Soit $A$ une partie définissable indécomposable d'un symétron $\omega$-stable $S$. Alors $s(A,x)$ est une partie définissable indécomposable pour tout $x\in S$.
\end{Lemme}
\begin{proof}
Soit $x\in S$; supposons par l'absurde que $s(A,x)$ n'est pas indécomposable. Il existe donc une famille finie d'ensembles convexes définissables deux-à-deux disjoints $A_1, ..., A_n$, telle que $s(A,x)\subseteq A_1\cup ...\cup A_n$ et $s(A,x)\cap A_i\neq \varnothing$ pour tout $i\in \{1, ..., n\}$, avec $n\geq 2$. Or, $s_x$ est une bijection définissable qui préserve la convexité. Par conséquent, $s(s(A,x),x)=A\subseteq s(A_1,x)\cup ...\cup s(A_n,x)$ et $A\cap s(A_i,x)\neq \varnothing$ pour tout $i\in \{1, ..., n\}$, contradiction.
\end{proof}
\begin{proposition}\label{indécomposbale} (à comparer avec \cite[Théorème 2.7]{Poi1})
Soit $A$ une partie définissable d'un symétron $\omega$-stable. Alors $A=X_1\cup...\cup X_n$, où les $X_i$ sont des parties définissables indécomposables deux-à-deux disjointes.
\end{proposition}
\begin{proof}
Si $A$ est indécomposable, on conclut. Sinon, il existe des parties convexes définissables deux-à-deux disjointes $X_1,..., X_n$ (avec $n\geq 2$) telles que \[A=(A\cap X_1)\cup...\cup(A\cap X_n),\]
où de plus, $A\cap X_i\neq \varnothing $, pour tout $1\leq i \leq n$.
Si par exemple $A\cap X_1$ n'est pas indécomposable, alors il existe des parties convexes définissables deux-à-deux disjointes $X_{1,1},...X_{1,k}$ (avec $k\geq 2$) avec :
\[A\cap X_1=(A\cap (X_1\cap X_{1,1}))\cup...\cup (A\cap (X_1\cap X_{1,k}))\]
où de plus, $(A\cap X_1)\cap X_{1,j}\neq \varnothing $, pour tout $1\leq j \leq k$.
Mais $X_1\cap X_{1,1}$ est une partie convexe définissable strictement contenue dans $X_1$ (sinon, $X_1\subseteq X_{1,1}$, et $X_{1,1}\cap X_{1,2}\neq \varnothing$, contradiction). Par la condition de chaîne descendante sur les parties convexes définissables, on construit un arbre à branchement fini et branches finies. Par le lemme de König, cet arbre est lui-même fini.
\end{proof}
\begin{Remarque}
Une partie définissable convexe connexe est indécomposable.
\end{Remarque}
Considérons maintenant une partie définissable $X$ d'un symétron $\omega$-stable $S$ et un sous-symétron (définissable) $S'$. On dit que $X$ \textit{engendre elliptiquement } $S'$ s'il existe un entier $n$ fixé, tel que tout élément de $S'$ est obtenu par au plus $n$ applications successives de prises de symétrie et de prises de milieu à partir d'éléments de $X$.
\begin{theoreme}\label{indécomposable symétron}(à comparer avec \cite[Théorème 14]{Poi4}) Soit $S$ un symétron de RM fini. Soit $(A_i : i\in I)$ une famille de parties définissables indécomposables telles que $\bigcap_{i\in I} A_i\neq \varnothing$. Alors $\bigcup_{i\in I} A_i$ engendre un symétron définissable connexe qui est en fait elliptiquement engendré par un nombre fini de $A_i$.
\end{theoreme}
\begin{proof}
Pour chaque suite $s=(i_1,...,i_n)$ et pour chaque $n$, on note $X_s=\bigcup_{i\in s} A_i$.
Soit $x\in \bigcap_{i\in I} A_i$. On considère $X$ de rang maximal parmi les parties définissables elliptiquement engendrées à partir d'un $X_s$ (on utilise ici que le rang de Morley de $S$ est fini). Soit maintenant $p$ un type générique de $X$. 
\\
Pour chaque $A_i$, l'ensemble $s_{A_i}(p)$  ne contient qu'un nombre fini de types $\{p_1,..., p_n\}$ de même rang que $p$ (on remarque que $s_x(p)$ appartient à cette famille finie). En effet, pour chaque $a_i\in A_i$, l'ensemble $s_{a_i}(X)$ est une partie définissable elliptiquement engendrée à partir d'un certain $X_s$ telle que $RM(s_{a_i}(X))=RM(X)$. Si $s_{A_i}(p)$ était infini, alors $s_{A_i}(X)$ serait une partie définissable elliptiquement engendrée par notre famille qui serait de rang strictement supérieur à celui de $X$, contradiction.
\\
On peut donc décomposer $A_i$ en un nombre fini de parties définissables deux-à-deux disjointes $A_{i,1}, ..., A_{i,n}$, où $A_{i,j}=A_i\cap Sym(p\longleftrightarrow p_j)$, $1\leq j \leq n$. Or, $A_i$ est indécomposable et donc $A_i=A_i\cap Sym(p\longleftrightarrow p_j)$ pour un certain $j$ avec $1\leq j \leq k$. Chaque $A_i$ est donc inclus dans le convexe définissable $S'=Sym(p\longleftrightarrow s_x(p))$. On a donc :
\[ RM(p)\leq RM(X)\leq RM(S')\leq RM(p).\] 
Finalement, $RM(X)=RM(S')$. Comme tous les $A_i$ sont dans la même composante connexe (par indécomposabilité, les $A_i$ sont contenus dans $Sym(p'\longleftrightarrow s_x(p'))$ pour $p'$ un type générique de $S'$), on peut supposer que $S'$ est connexe. Finalement, tout élément de $S'$ est au milieu de deux points génériques et donc de deux points de $X$; par conséquent, $S'=m(X,X)$ et $S'$ est elliptiquement engendré par les $A_i$, $i\in s$.
\end{proof}

\section{Une analyse de Lascar pour les boucles et les symétrons}
\subsection{Connexité}
On peut tout d'abord établir une correspondance entre certaines parties convexes définissables et les sous K-boucles définissables.
\begin{proposition} \label{convexe et boucle} Soit $B$ une K-boucle uniquement 2-divisible $\omega$-stable. Une partie définissable $X\subseteq B$  est une sous K-boucle ssi elle est convexe et contient 0. De plus, une sous K-boucle définissable est uniquement 2-divisible.
\end{proposition}
\begin{proof}
Le sens direct est clair. Supposons que $X$ est une partie définissable convexe (et donc aussi close par prise de milieu) contenant 0. On remarque que $X$ est stable par inversion : pour $x\in X$, on a $s(x,0)=0+(-x+0)=-x$.  D'après la Proposition \ref{sous-boucle}, il suffit de montrer que $X$ est stable par addition. Remarquons tout d'abord que $X$ est stable par passage au double et par extraction des moitiés : $s(0,x)=x\cdot 2$ et $s(0,x/2)=x$.
\\Soit $x,y\in X$, on considère $(x+y)\cdot 2$. Or, d'après le Fait \ref{boucle général}, on a $(x+y)\cdot 2=x+(y\cdot 2+x)=s((-y)\cdot 2,x)$, mais $(-y)\cdot 2\in X$, et donc $(x+y)\cdot 2\in X$. Finalement, par extraction des moitiés, on obtient que $(x+y)\in X$. 
\\
Soit $B_1$ une sous K-boucle définissable; c'est une partie convexe définissable donc close par prise de milieu. Pour tout $b\in B_1$, l'élément $b/2$ est le milieu de $b$ et $0$; on a bien $b/2\in B_1$.
\end{proof}
\begin{corollaire}
On peut associer à chaque partie convexe définissable $X$ une sous-boucle définissable de même rang (à savoir $s(X,x/2)$ pour $x\in X$).
\end{corollaire}
\begin{Remarque}
L'image d'une partie convexe définissable par un automorphisme définissable est une partie convexe définissable. En particulier, l'image d'une partie indécomposable par un automorphisme définissable reste indécomposable. 
\end{Remarque}
Nous pouvons désormais clarifier la notion de connexité dans notre contexte.
\begin{definition}
 Une K-boucle uniquement 2-divisible $\omega$-stable est b-connexe si elle ne contient pas de sous-boucle définissable propre de même rang.
\end{definition}
\begin{proposition}
Soit $B$ une K-boucle uniquement 2-divisible $\omega$-stable. La boucle $B$ est b-connexe ssi elle est connexe (en tant que symétron) ssi son degré de Morley est égal à un.
\end{proposition}
\begin{proof}
Supposons que le degré de Morley de $B$ est égal à un. Si la boucle contenait un sous-symétron $X$ propre définissable de même rang alors $B\setminus X$ serait également une partie générique définissable disjointe de $X$ (car $X$ s'injecte dans $B\setminus X$), contradiction.
\\
Si $B$ est connexe alors $B$ est b-connexe : une sous-boucle propre définissable de même rang est en particulier un sous-symétron propre définissable de même rang.
\\
Si $B$ est $b$-connexe alors le degré de $B$ est égal à un. Sinon, il y aurait au moins deux types génériques distincts.  Soient $p_1,..,p_n$ les types génériques, pour $1\leq i \leq n$, les parties de la forme $Sym(p_1\longleftrightarrow p_i)$ sont des parties convexes définissables qui contiennent un unique générique et qui sont deux-à-deux disjointes. Une de ces parties contient 0 et forme donc une sous-boucle définissable propre de même rang, contradiction.
\end{proof}
\begin{definition}
Soit $B$ une K-boucle uniquement 2-divisible $\omega$-stable. On dira que la sous-boucle connexe de même rang que $B$ est la composante connexe (principale) de $B$ et on la notera $B^{\circ}$.
\end{definition}
Si la boucle $B$ est en plus de RM fini, on obtient une version particulièrement commode du théorème des indécomposables.
\begin{corollaire}\label{indécomposable boucle}
Soit $B$ une K-boucle uniquement 2-divisible de RM fini et soit $(A_i : i\in I)$ une famille de parties définissables indécomposables contenant $0$. Alors $\bigcup_{i\in I} A_i$ engendre une sous K-boucle connexe définissable qui est en fait elliptiquement engendrée par un nombre fini de $A_i$.
\end{corollaire}
\subsection{Analyse de Lascar}
Nous pouvons dès lors tenter d'adapter l'analyse de Lascar aux K-boucles uniquement 2-divisibles de RM fini. Le principal obstable réside dans la construction d'une sous-boucle normale à partir d'ensembles indécomposables. En effet, dans une boucle, il n'y a pas forcément de groupe définissable d'automorphismes intérieurs.
\begin{definition}
On dit qu'une boucle $B$ est automorphique si $I(B)$, le groupe des applications intérieures (Définition \ref{groupe application intérieure}), est composé d'automorphismes. 
\end{definition}
On rappelle qu'une sous-boucle est normale ssi elle est stabilisée par le groupe $I(B)$ des applications intérieures.
\\
La propriété d'être une K-boucle uniquement 2-divisible automorphique passe au quotient :
\begin{Lemme}
Soit $B$ une K-boucle uniquement 2-divisible automorphique $\omega$-stable et soit $B_1$ une sous-boucle définissable normale. Alors la boucle quotient $B/B_1$ est une K-boucle uniquement 2-divisible automorphique.
\end{Lemme}
\begin{proof}
La boucle quotient $B/B_1$ est une K-boucle car les propriétés équationnelles sont préservées par passage au quotient. De plus, elle est bien 2-divisible car  $\overline{a}=\overline{a'+a'}=\overline{a'}+\overline{a'}$ pour $a=a'+a'$. On remarque que si $\overline{a}\cdot 2=\overline{0}$ alors $a\cdot 2\in B_1$. Mais $B_1$ est une partie convexe définissable; c'est donc un ensemble clos par prise de moitié; on a donc $a\in B_1$. La K-boucle quotient $B/B_1$ ne contient pas d'élément d'ordre 2 et par le Fait \ref{boucle général}, on en déduit que l'application qui à $\overline{a}$ associe $\overline{a}\cdot 2$ est injective.
\\
\\
Pour ce qui est du fait d'être automorphique, on utilise \cite[Théorème 2.2]{Bruck}.
\end{proof}
\begin{theoreme}\label{Lascar} (à comparer avec \cite[Théorème 2.13]{Poi1})
Soit $B$ une K-boucle uniquement 2-divisible de RM fini ($\omega$-saturée) qui est de plus automorphique. Alors il existe des formules fortement minimales $\phi_1, ..., \phi_n$ telles que dans toute extension élémentaire propre, l'ensemble des réalisations d'au moins une de ces formules augmente.
\end{theoreme}
\begin{proof}
 Soit $A$ une partie définie par une formule fortement minimale. D'après la Proposition \ref{indécomposbale}, il existe une partie définissable $A'$, de même rang et même degré qui est indécomposable. Or, pour $x\in A'$, $s(A',x/2)=A_1$ est une partie définissable fortement minimale indécomposable qui contient $0$.
\\
On considère la famille $(f(A_1) : f\in I(B))$. Remarquons que pour tout $f\in I(B)$, $f(s(a,a'))=s(f(a),f(a'))$ et $f(m(a,a'))=m(f(a),f(a'))$ pour $a, a'\in A_1$. Appliquons le Théorème \ref{indécomposable symétron} à la famille $(f(A_1) : f\in I(B))$ : on obtient une sous K-boucle définissable connexe normale $B_1$ elliptiquement engendrée à partir d'une partie définissable fortement minimale (la sous-boucle ainsi obtenue est bien quasi-fortement minimale). 
\\
Si $B_1=B^{\circ}$, on arrête; sinon, on considère la K-boucle interprétable $B^{\circ}/B_1$ (par connexité, $B_1\subseteq B^{\circ}$) et on repète le même procédé.
\\
Nous obtenons une suite de sous K-boucles définissables connexes : $B_1\triangleleft ...\triangleleft B_n=B^{\circ}$. En passant à une extension élémentaire propre, l'ensemble des réalisations d'une des formules fortement minimales déterminant les quotients doit augmenter.
\end{proof}
\begin{Remarque}
On va voir que les formules fortement minimales obtenues dans le théorème précédent peuvent être prises à paramètres dans le modèle premier de $Th(B)$. L'hypothèse sur la saturation n'est donc pas nécessaire.
\end{Remarque}
\begin{definition}
Soit $T$ une théorie $\omega$-stable et soient $p,q$ deux types. On dit que $p$ et $q$ sont orthogonaux si pour tout ensemble $A$ contenant une base de $p$ et de $q$, toute réalisation $\overline{a}$ d'une extension non-déviante de $p$ à $A$ est indépendante de toute réalisation $\overline{b}$ d'une extension non-déviante de $q$ à $A$.
\end{definition}
\begin{definition}
Soit $T$ une théorie $\omega$-stable. On dit qu'elle est fini-dimensionnelle s'il existe une famille finie de types stationnaires telle que toute type non-algébrique est non-orthogonal à un des types de la famille.
\end{definition}
\begin{fait}\label{modèle premier} \cite[Proposition 3. Sec.2. Chap.8]{Las}
Soit $A\subset \mathfrak{C}$, alors il existe un modèle $\mathfrak{M}$ absolument premier sur $A$, \textit{i.e.}, il existe une énumération de $M=\{ m_{\alpha} : \alpha<\lambda \}$ (où $\lambda$ est un ordinal) telle que pour tout $\alpha <\lambda$, $tp(m_{\alpha}/A\cup \{m_{\beta} : \beta<\alpha \})$ est isolé. En particulier, $\mathfrak{M}$ est premier sur $A$.
\end{fait}
\begin{corollaire}
   Soit $B$ une K-boucle uniquement 2-divisible de RM fini. Si $B$ est automorphique, alors $Th(B)$ est fini-dimensionnelle. De plus, les formules fortement minimales qui déterminent les dimensions peuvent être prises à paramètres dans le modèle premier.
\end{corollaire}
\begin{proof}
Soit $\mathfrak{B}$ un modèle monstre de $Th(B)$ et soit $A\subset \mathfrak{B}$. Soient $\phi_1, ..., \phi_n$ les formules fortement minimales obtenues dans le Théorème \ref{Lascar}. On considère un type non-algébrique $p$ sur $A$. 
\\
Prenons un modèle $B$ suffisamment saturé contenant $A$ et les paramètres des $\phi_i$, pour $i\in \{1, ..., n \}$. Soit $p'\in S_n(B)$ une extension non-déviante de $p$ et $\overline{a}\models p'$. D'après le Fait \ref{modèle premier}, il existe un modèle $B'$ absolument premier sur $B\cup \{\overline{a} \}$. Or, d'après le Théorème \ref{Lascar}, on peut supposer qu'il existe $\overline{b}\in \phi_1[B']-\phi_1[B]$. Mais $tp(\overline{b}/B\cup{\overline{a}})$ est isolé et donc $\overline{a}\nind_{B} \overline{b}$. Par conséquent, le type $p$ n'est pas orthogonal au type fortement minimal déterminé par $\phi_1$.
\\
Une théorie fini-dimensionnelle n'a pas la propriété du recouvrement fini et donc une formule minimale est fortement minimale (voir \cite[Théorème 2.13]{Poi1} ou bien \cite[Proposition 4.7.4 et Lemme 4.7.5]{Wagl}). Il suffit de procéder comme dans la démonstration du Théorème \ref{Lascar} en considérant le modèle premier de $Th(B)$.
\end{proof}
\begin{fait}\cite[Corollaire 2.14]{Poi1} ou \cite[Proposition 4.7.10]{Wagl}
Dans une théorie fini-dimensionnelle où les dimensions sont portées par des types fortement minimaux, le rang de Morley est fini et il égal au rang de Lascar. De plus, pour tout entier $n$, l'ensemble $\{a : RM(\phi[x,a])=n \}$ est définissable.
\end{fait}
\begin{corollaire} 
  Dans une K-boucle uniquement 2-divisible de RM fini qui est de plus automorphique, le rang de Lascar et le rang de Morley coïncident. C'est une K-boucle rangée au sens de Borovik.
\end{corollaire}
\begin{Remarque}
Les symétrons et les K-boucles qui apparaissent dans la classification des groupes de RM fini sont en général inteprétables dans le groupe ambiant et sont donc rangés. Néanmoins, les résultats précédents conservent leur pertinence lorsque l'on s'intéresse aux symétrons et aux K-boucles en tant que tels sans référence à un groupe ambiant de rang de Morley fini.
\end{Remarque}
B. Poizat a reformulé les résultats précédents dans le seul langage des symétrons \cite{Poi5}. 
En guise de conclusion, nous indiquons quelques directions possibles de recherche :
\begin{enumerate}
    \item Est-il possible de lever l'hypothèse concernant le caractère automorphique des K-boucles dans l'analyse de Lascar ?
    \item Peut-on donner une classification des symétrons \textit{algébriques}, i.e l'ensemble de base est une variété algébrique et l'opération de symétron est un morphisme ? En particulier, un symétron définissable dans un pur corps algébriquement clos est-il un symétron algébrique ?
    \item En lien un avec la question précédente, un symétron donné génériquement dans une théorie $\omega$-stable donne-t-il lieu à un symétron définissable ? 
\end{enumerate}
Dans \cite{Poi6}, B. Poizat dresse un panaorama exhaustif de l'état de nos connaissances sur les symétrons ainsi que des questions qui restent à élucider.

\end{document}